\documentclass[12pt]{amsart}
\usepackage{amsmath,amsthm,amsfonts,amssymb}
\begin{document} 
\newcommand{\B}{{\mathbb B}}
\newcommand{\C}{{\mathbb C}}
\newcommand{\N}{{\mathbb N}}
\renewcommand{\O}{{\mathcal O}}
\newcommand{\A}{{\mathbb A}}
\newcommand{\Q}{{\mathbb Q}}
\newcommand{\Z}{{\mathbb Z}}
\renewcommand{\P}{{\mathbb P}}
\newcommand{\R}{{\mathbb R}}
\newcommand{\rc}{\subset}
\newcommand{\rank}{\mathop{rank}}
\newcommand{\trace}{\mathop{tr}}
\newcommand{\tensor}{\otimes}
\newcommand{\dimc}{\mathop{dim}_{\C}}
\newcommand{\Lie}{\mathop{Lie}}
\newcommand{\Auto}{\mathop{{\rm Aut}_{\mathcal O}}}
\newcommand{\alg}[1]{{\mathbf #1}}
\newtheorem*{definition}{Definition}
\newtheorem*{maintheorem}{Main theorem}
\newtheorem*{claim}{Claim}
\newtheorem{corollary}{Corollary}
\newtheorem*{Conjecture}{Conjecture}
\newtheorem*{SpecAss}{Special Assumptions}
\newtheorem{example}{Example}
\newtheorem*{remark}{Remark}
\newtheorem*{observation}{Observation}
\newtheorem*{fact}{Fact}
\newtheorem*{remarks}{Remarks}
\newtheorem{lemma}{Lemma}
\newtheorem{proposition}{Proposition}
\newtheorem{theorem}{Theorem}
\title[Brody Curves]{%
A projective manifold where entire and Brody curves
behave very differently.
}
\author {J\"org Winkelmann}
\begin{abstract}
We give an example of a projective manifold with dense entire
curves such that every Brody curve is degenerate.
\end{abstract}
%
\address{%
J\"org Winkelmann \\
 Institut Elie Cartan (Math\'ematiques)\\
 Universit\'e Henri Poincar\'e Nancy 1\\
 B.P. 239\\
 F-54506 Vand\oe uvre-les-Nancy Cedex\\
 France
}
\email{jwinkel@member.ams.org\newline\indent{\itshape Webpage: }%
http://www.math.unibas.ch/\~{ }winkel/
}
\thanks{
{\em Acknowledgement.}
We are grateful for useful discussions with Serge Cantat and
Frederic Campana.
}
\maketitle
\section{Introduction}
Let $X$ be a complex manifold.
An entire curve is a non-constant holomorphic map from $\C$ to $X$.
Assume $X$ to be endowed with a hermitian metric.
Then an entire curve $f:\C\to X$ is called a ``Brody curve''
iff its derivative is bounded with respect to the euclidean metric
on $\C$ and the given hermitian metric on $X$. If $X$ is compact,
the notion of a ``Brody curve'' is independent of the choice of
the metric.

The famous result of Brody (\cite{B}) implies:
{\em A compact complex manifold $X$ admits an entire curve if and
only if it admits a Brody curve.}

This is an important result which, for example, allows an easy
characterization of those submanifolds of abelian varieties which
are hyperbolic in the sense of S.~Kobayashi.

Given the above result, it is natural to ask the following question:
{\em ``Let $X$ be a compact complex manifold and $p\in X$.
Assume that there exists an entire curve $f:\C\to X$ with
$p\in f(\C)$. Does this implies that there is a Brody
curve $f$ with $p\in f(\C)$?''}

In \cite{W} we gave an example of a non-compact manifold
(actually a domain in an abelian variety) where the behaviour
of entire curves and Brody's curves differ. 

However this was not really a negative answer, since the question
really concerns compact manifold.

Now we are able to give an example of a compact complex manifold
where the above formulated question has a negative answer.

We prove the following theorem:
\begin{maintheorem}
There exists a projective manifold $X$ with a hypersurface $Z$
such that for every point $x\in X$ there
exists an entire curve $f:\C\to X$ with $\overline{f(\C)}=X$
and $x\in f(\C)$,
but $Z$ contains the image of every Brody curve.
\end{maintheorem}

\section{Proof of the Main theorem}
The main theorem is a consequence of the more specific theorem below.
\begin{theorem}
There exists an abelian threefold $A$ with a smooth curve $C$
such that the smooth projective variety $\hat A$ obtained by
blowing up $A$ along $C$ has the following properties:
\begin{enumerate}
\item
For every point $p\in \hat A$ there exists a non-constant
entire curve $\gamma:\C\to\hat A$ with $p\in\gamma(\C)$
and $\overline{\gamma(\C)}=\hat A$.
\item
For a given point $p\in \hat A$ there exists a non-constant
Brody curve $\gamma:\C\to\hat A$ with $p\in\gamma(\C)$
if and only if $p$ is contained in the exceptional divisor
$E=\pi^{-1}(C)$ of the blow-up $\pi:\hat A\to A$.
\end{enumerate}
\end{theorem}

The key idea is that the hermitian metric will explode in some
directions due to the blow up and that this will create an
obstruction against lifting Brody curves. As a consequence
there will be no Brody curves outside the exceptional divisor.
To realize this idea it is necessary to ensure that the center
of the blow-up intersects the closure of the image of each
brody curve. To achieve this it will be necessary to blow up
a center of positive dimension.
Furthermore, since the center of a blow-up has real codimension
at least four, it is necessary to choose the abelian variety
in such a way that for every Brody curve the closure is at
least real $4$-dimensional.

\begin{proof}
$(1)$. For every $x\in A$ and $v\in T_xA$ there is an affine-linear
curve $\gamma:\C\to A$ with $\gamma(0)=x$ and $\gamma'(0)=v$.
Recall that $\pi:\hat A\to A$ is an isomorphism outside $C$ and that
each point of $x\in C$ is replaced by $\P(T_x/T_xC)$. Observe further
that each entire curve $\gamma:\C\to A$ lifts to $\hat A$
unless $\gamma(\C)\subset C$. Combined, these facts yield
statement $(1)$.

$(2)$.
We have $\pi^{-1}(x)\simeq\P_1$ for every $x\in C$.
This implies that there is a Brody curve through every point
in $E=\pi^{-1}(C)$.
Conversely, let $\hat f:\C\to\hat A$ be a Brody curve.
We will see that $\hat f(\C)\subset E=\pi^{-1}(C)$ if
we choose $A$ and $C$ according to prop.~\ref{prop-ex}
below.
Now $\hat f$ being a Brody curve implies
that $f=\pi\circ \hat f:\C\to A$ is a Brody curve or constant.
Let us assume that $f$ is not constant.
If $A=\C^3/\Gamma$, then $f$ lifts to an affine-linear map
$F:\C\to \C^3$. $f(\C)$ is thus the orbit of a complex
one-parameter subgroup $P$ of $A$. Let $H$ denote the
(real) closure of $P$ in $A$.
Thanks to prop.~\ref{prop-ex} we may assume that $H$ and $C$
intersect transversally in some point $p$.
But now we arrive at a contradiction because according
to prop.~\ref{brody-blowup} under these circumstances $f:\C\to A$
can not be induced by a Brody curve $\hat f:\C\to\hat A$.
Thus $f=\pi\circ\hat f$ must be constant. Since $\hat f$
is non-constant and $\pi$ is an isomorphism outside of $C$,
it follows that $\hat f(\C)\subset E=\pi^{-1}(C)$.
\end{proof}
\section{Brody curves}
We recall some basic facts on Brody curves.

Let $X$ be a complex manifold endowed with some hermitian metric.
Then an {\em entire curve} is a non-constant holomorphic map
from $\C$ to $X$ and a {\em Brody curve} is a non-constant
holomorphic map $f:\C\to X$ for which the derivative $f'$ is
bounded (with respect to the euclidean metric on $\C$ and the
given hermitian metric on $X$).

If $X$ is compact, the notion of a ``Brody curve'' is independent
of the choice of the hermitian metric.

If $\phi:X\to Y$ is a holomorphic map between compact complex
manifolds and $f:\C\to X$ is a Brody curve, then
$\phi\circ f:\C\to Y$ is a Brody curve, too.
(But not necessarily conversely.)

If $X=\C^g/\Gamma$ is a compact complex torus (e.g.~an abelian
variety), then an entire curve $f:\C\to X$ is a Brody curve
if and only if it lifts to an affine linear map
$\tilde f:\C\to \C^g$.

\section{Local model of blow-up}
The idea we use is:
If we blow up something, the hermitian metric will explode
somewhere.
We will now make this precise.
\begin{proposition}\label{metric-blowup}
Let $A$ be a three-dimension complex manifold, $C$ a smooth curve,
$\pi:\hat A\to A$ the corresponding blow-up with center $C$,
$p\in C$ and $L_n$ a sequence of curves converging to a
curve $L_0$ such that
\begin{enumerate}
\item
$L_0$ intersects $C$ transversally in $p$.
\item
The intersection $L_n\cap C$ is empty for all $n\ne 0$.
\end{enumerate}
Furthermore assume $A$ and $\hat A$ endowed with hermitian metrics.

Then there exists  sequences $p_n\in L_n$ and $v_n\in T_{p_n}(L_n)$
such that $\lim p_n=p$ and
\[
\limsup \frac { ||\pi^{-1}(v_n)||_{\hat A}}
{||v_n||_{A}}=\infty
\]
(note that $L_n\subset A\setminus C$ and that $\pi$ is an isomorphism
on $A\setminus C$.)
\end{proposition}
\begin{proof}
Let $\hat p$ denote the point in $\pi^{-1}(p)$ which points in the
direction of $L$ (using the isomorphism between $\pi^{-1}(p)$
and the projectivization of the normal tangent space $T_pA/T_pC$.)

We fix local holomorphic coordinates on $A$ and $\hat A$ around
$p$ resp.~$\hat p$ such that the defining equations for $C$ and $L_0$
become as simply as possible. Doing this we get 
local holomorphic coordinates such that
\[
C=\{(z_1,z_2,z_3):z_1=z_2=0\},
\]
\[
L_0=\{(z_1,z_2,z_3):z_2=z_3=0\},
\]
Now the projection:
\[
\pi(x_1,x_2,x_3)=(x_1x_2,x_2,x_3)
\]
Since $\lim L_n=L_0$, the curves $L_n$ can be parametrized
as
\[
L_n=\{\gamma_n(t)=(t,\alpha_n(t),\beta_n(t)\}
\]
where $t$ runs through an appropriate small neighbourhood of $0$
and where $\alpha_n$, $\beta_n$ are sequences of holomorphic functions
converging uniformly to the constant function zero on this
small neighbourhood.

Since all the calculations happen in some small neighbourhood of $p$
resp.~$\hat p$, we may replace the given hermitian metrics by the
euclidean metric with respect to our coordinate systems.

Our next step is to define the auxiliary function
\[
\phi_n(t)=t+\alpha_n(t)\alpha_n'(t)
\]
We observe that $\phi_n$ converges to the identity map $\phi(t)=t$.
Therefore the theorem of Rouche allows us to choose a sequence
$s_n$ with $\lim_n s_n=0$ and $\phi_n(s_n)=0$ for all $n$.

We claim:
{\em $\alpha_n(s_n)\ne 0$}. Indeed, assume $\alpha_n(s_n)=0$.
Then 
\[
\alpha_n(s_n)\alpha_n'(s_n)=0
\]
 and consequently
\[
0=\phi_n(s_n)=s_n+0\ \Rightarrow\ s_n=0
\]
and therefore
\[
0=\alpha_n(s_n)=\alpha_n(0).
\]
But $\alpha_n(0)=0$ is impossible, because $L_n\cap C$ is empty.
Thus the assumption $\alpha_n(s_n)=0$ leads to a contradiction,
i.e. $\alpha_n(s_n)$ must be non-zero.

Hence we may divide by $\alpha_n(s_n)$ and thereby deduce that
$\phi_n(s_n)=0$ implies $\alpha_n'(s_n)=-s_n/\alpha_n(s_n)$.
If $\hat\gamma(t)$ denotes the point in $\hat A$ lying above
$\gamma_n(t)\in A\setminus C$, we obtain
\[
\hat\gamma_n(s_n)=
\left( \frac{s_n}{\alpha_n(s_n)}, \alpha_n(s_n), \beta_n(s_n) \right)
=
\left( -\alpha_n'(s_n), \alpha_n(s_n), \beta_n(s_n) \right)
\]
which converges to $(0,0,0)=\hat p\in\hat A$ if $n$ goes to infinity.

Now
\[
\gamma_n'(s_n)=(1,\alpha_n'(s_n),\beta_n'(s_n))
\ \ \Rightarrow\ \ 
\lim_n||\gamma_n'(s_n)||=1
\]
while
\[
\hat\gamma_n'(s_n)=\left(\frac{1-\alpha_n'(s_n)}{(\alpha_n(s_n))^2},
\alpha_n'(s_n),\beta_n'(s_n)\right)
\ \ \Rightarrow\ \ 
\lim_n||\hat\gamma_n'(s_n)||=+\infty
\]
\end{proof}
\section{Brody curves and blow ups}
\begin{proposition}\label{brody-blowup}
Let $A$ be an abelian variety, $C$ a submanifold 
containing $e_A$ and $\pi:\hat A\to A$
the blow-up with center $C$. 

Let $\gamma:\C\to A$ be a Brody curve with closure
$B=\overline{\gamma(\C)}$ such that 
$\gamma(0)=e_A\in C$.
Assume moreover that
\begin{enumerate}
\item
$v=\gamma'(0)\not\in T_eC$ and
\item
$T_eB\not\subset T_eC\oplus\left<v\right>_{\C}$
\end{enumerate}

Then there does not exist a Brody curve $\hat\gamma:\C\to\hat A$
with $\gamma=\pi\circ\hat\gamma$.
\end{proposition}
\begin{proof}
Let $\tau:\C^3\to A$ be the universal covering.
Since $\gamma$ is a Brody curve with $\gamma(0)=e$, 
it lifts to a linear
map $\tilde\gamma:\C\to \C^3$. 

Recall that $B$ is a real subtorus. Let $V$ be the Lie algebra
of $B$. We may regard $V$ as the connected component of $\pi^{-1}(B)$
which contains $0$. Since $\gamma(\C)$ is dense in $B$, we can find
a finitely generated subgroup $\Lambda_0\subset\pi^{-1}(\gamma(\C))$
which is dense in $V$. Define $\Lambda=\pi(\Lambda_0)$.
Since $T_pB\not\subset T_pC+\left<\gamma'(0)\right>_{\C}$ 
we can choose a sequence
$\lambda_n\in\Lambda$ such that $\lim\lambda_n=e$ and 
$C$ does not intersect the translate $\lambda_n+\gamma(\Delta)$
(where $\Delta$ is a sufficiently small disk containing $0$)
in some neighbourhood of $p$.

Now we can invoke prop.~\ref{metric-blowup}
taking $X$ to be some small open neighbourhood of $e$ in $A$,
$p=e_A$, $L_0=X\cap\gamma(\Delta)$
and $L_n=X\cap(\lambda_n+\gamma(\Delta))$.

Hence
\begin{equation}\label{eq}
\sup_{t\in C}\frac{||\hat\gamma'(t)||_{\hat A}}%
{||\gamma'(t)||_A}=+\infty
\end{equation}
where $\hat\gamma:\C\to\hat A$ is the natural lift of $\gamma$.

Since $\gamma$ is induced by an affine-linear map, the norm
$||\gamma'(t)||$ is a positive constant and in particular bounded
from below by a number greater than zero.
Together with the above equation \ref{eq} this implies
that $\hat\gamma$ can not be Brody curve.
\end{proof}

\section{Excluding Real subtori of dimension three}
In this section we deduce the following statement:
\begin{proposition}\label{3-subtori}
There exists an abelian three-fold $A$ such that every real
subtorus of real dimension three is totally real in $A$.
\end{proposition}

We will prove this assertion by showing that every very general
abelian three-fold has this property, i.e. we demonstrate:
\begin{proposition}\label{prop-uni}
Let $U\to D$ be a locally complete family
of abelian varieties of dimension three.

Then there exists a countable family of nowhere dense
closed analytic subsets $Z_i\subset D$ such that every abelian
threefold $A$ corresponding to a point outside the
union $\cup_i Z_i$ has the property
``Every real
subtorus of real dimension three is totally real in $A$''
\end{proposition}

Before proving the proposition, we need some lemmata.

\begin{lemma}\label{ex-def}
Let $A=\C^3/\Gamma$ be an complex abelian $3$-fold,
$S$ a real subtorus of dimension three.

Then there is a joint deformation of $S\subset A$ over the
unit disc such that $A_t$ is an abelian variety for all $t$
and $S_t$ is totally real for all $t\ne 0$.
\end{lemma}
\begin{proof}
Let $\Lambda\subset\Gamma$ be the $\Z$-submodule 
corresponding to $S$. Since $A$ is an abelian variety,
$\C^3$ admits a hermitian form $H$ such that
$B=\Im H$ has integer values on $\Gamma\times\Gamma$.
Now $B$ is alternating and $3=\rank_\Z(\Lambda)$ is odd,
hence there is an element $v\in\Lambda$ for which
$B(v,\cdot)$ vanishes identically on $\Lambda$.
Let $\Lambda_\R$ resp.~$\Lambda_\C$
be the real resp.\ complex vector subspace of $\C^3$ generated
by $\Lambda$.
We may assume that $\Lambda_\R$ is not totally real. Then
$\dim_\C(\Lambda_\C)=2$ and $L=\Lambda_\R\cap i\Lambda_\R$ is a 
complex line.
Now we choose an element $w\in \Gamma$ such that
\begin{enumerate}
\item
$B(v,w)\ne 0$,
\item
$B(w,\cdot)$ does not vanish identically on $L$ and
\item
$w\not\in\Lambda_\C$.
\end{enumerate}
We define $\R$-linear self-maps $\phi_t$ of $\C^3$ as follows:
First we observe that $\C^3$ is the direct sum of $\R\cdot v$ and
$K=\{x:B(x,w)=0\}$.
Second we set
$\phi_t(v)=v+tw$ and $\phi_t(x)=x$ for all $x\in K$.
It is easy to check that $\phi_t$ is always bijective and
moreover an isometry for $B$.
Hence $\Gamma_t=\phi_t(\Gamma)$ is a lattice for which
the assertion $B(\Gamma_t,\Gamma_t)\subset\Z$ holds.
Thus $A_t=\C^3/\Gamma_t$ is an abelian variety.

Now let us look at $\phi_t(\Lambda_\R)$.
First we consider the real vector subspace 
$V=\Lambda_\R\oplus\R w$. Let $K=\{x:B(x,w)=0\}$.
Then $V=(V\cap K)\oplus\R v$. Now $\phi_t$ acts trivially on $K$
and $\phi_t(v)=v+tw\in V$. Hence $\phi_t$ stabilizes $V$.
We note that $\dim_\R(V)=4$ and $<V>_\C=\C^3$, because $w\not\in
\Lambda_\C$. Therefore $V$ contains a unique complex line,
which must be $L$. Since $\phi_t(\Lambda_\R)\subset V$,
we may deduce that for each $t$ either $\phi_t(\Lambda_\R)$
is totally real or contains $L$. 
Now, by the construction of $\phi_t$ it is clear that
\[
\phi_t(\Lambda_\R)\cap\phi_s(\Lambda_\R)=\Lambda_R\cap K
\]
for any $s\ne t$. Since $L\ne\Lambda_\R\cap K$ due to
condition $(2)$ for the choice of $w$, we may deduce
that $L\not\subset\phi_t(\Lambda_R)$ for $t\ne 0$.
As a consequence, $\phi_t(\Lambda_\R)$ is totally real
for $t\ne 0$.
\end{proof}
\begin{lemma}\label{ntr-ana}
Let $\pi:U\to D$ be a family of three-dimensional complex
abelian varieties, parametrized by $D$ which we assume
to be the unit ball in some $\C^N$.

Let $S_0$ be a real three-dimensional subtorus of
the abelian variety $U_0=\pi^{-1}(0)$

Then there is natural deformation $S_t$ of $S_0$ ($t\in D$)
such that 
\[
Z=\{t\in D: S_t \text{ is not totally real}\}
\]
is a closed complex analytic subset of $D$.
\end{lemma}
\begin{proof}
The family $U$ can be described as a quotient $\C^3\times D$
by a $\Z^{6}$-action which is given as
\[
(m_1,\ldots,m_3,n_1,\ldots,n_3):
(v;t)\mapsto \left(v+(m_1,m_2,m_3)+ \sum_i n_if_i(t);t\right)
\]
where $v\in \C^3$, $t\in D$ and where the $f_i$ are holomorphic maps
from $D$ to $\C^3$.

We may assume that $S_0$ is the subtorus fro which the corresponding
subgroup of $\Z$-rank $3$ is generated by $f_1(0)$, $f_2(0$ and
$f_3(0)$. Then $S_t$ corresponds to the subgroup generated
by the $f_i(t)$ and $S_t$ is totally real if and only if this
group spans $\C^3$ as a complex vector space. Therefore the set of
all $t\in D$ for which $S_t$ fails to be totally real is the
zero locus of $\det(f_1(t),f_(t),f_3(t))$ and thus a closed
complex analytic set.
\end{proof}
Now we can prove the proposition
\begin{proof}
There are only countably many different real subtori of real
dimension three for a given abelian $3$-fold $A$, each corresponding
to a $\Z$-submodule of rank three of $\Z^6=H_1(A,\Z)$.

Inside the family $U\to D$ there are canonical isomorphisms
\[
H^1(U_0,\Z)\simeq H^1(U_t,\Z)
\]
 which we may therefore identify.

Now for each fixed $\Z$-submodule of rank three the set of all $t\in D$
for which the corresponding subtorus $S_t$ fails to be totally real
is a closed analytic subset (lemma~\ref{ntr-ana}) which is not
all of $D$ (lemma~\ref{ex-def}).
This proves the proposition \ref{prop-uni} and thereby
prop.~\ref{3-subtori}.
\end{proof}
\begin{remark}
$1.)$ Since every subtorus of dimension smaller than three can be embedded
into a subtorus of dimension three, the property
{\em ``All real subtori of real dimension three
are totally real''} is equivalent to the property
{\em ``All real subtori of real dimension up to three
are totally real''}.
$2.)$ An abelian threefold $A$ is a {\em simple} abelian variety
iff it contains no elliptic curve. The latter property is equivalent
to the statement {\em ``All real subtori of real dimension up to two
are totally real''}. Hence the property {\em ``All real subtori of 
real dimension three
are totally real''} implies that the abelian $3$-fold under discussion
is simple.
\end{remark}
\section{Dealing with real subtori of dimension four}
The main goal of this section is to to verify that we can a choose
a curve $C$ in a $3$-dimensional abelian variety $A$ such that $C$
intersects the closure of every translate of every real subtorus
of real dimension four.

\begin{lemma}\label{bertini}
Let $A$ be an abelian threefold, $x\in A$ and let $L$
be a complex line in $T_xA$.

Then there exist smooth curves $C\subset A$ with $x\in C$
such that $T_xC$ is arbitrarily close to $L$.
\end{lemma}
\begin{proof}
We construct curves by embedding $A$ into a projective space
and taking the intersection of $A$ with linear subspaces of
codimension two containing $x$. Then the statement
follows from Bertini's theorem.
\end{proof}
\begin{lemma}\label{union}
Let $A$ be an abelian threefold with smooth curves $C$ and $C'$.

Then there is a dense open subset $U\subset A$ such that 
$C\cup\lambda_t^*C'$ is smooth for $t\in U$ where
$\lambda_t$ denotes translation by $t$.
\end{lemma}
\begin{proof}
$C\cup\lambda_t^*C'$ is smooth iff $C$ and $\lambda_t^*C'$ are
disjoint. Hence $U=A\setminus\{x-y:x\in C,y\in C'\}$.
\end{proof}

\begin{lemma}\label{linalg}
Let $V$ be a complex three-dimensional vector space equipped with
a hermitian inner product $H$ and let $W$ be a real four-dimensional
real subspace. Then there exists a complex line $L\subset V$
such that the angle between $W$ and $L$ is at least $\pi/4$.
i.e.,
\[
|<v,w>|\le\cos(\pi/4)||v||\cdot||w||
\]
for all $v\in L,w\in W$.
\end{lemma}
{\em Remark.}
If $H(\ ,\ )$ is an hermitian inner product, its real part is
the associated euclidean inner product and thus the angle between
two vectors $v$ and $w$ is the number $\phi\in[0,\pi/2]$ for which
$\cos\phi=\Re H(v,w)$.
\begin{proof}
We may choose vectors $A,B,C$ such that $(A,iA,B,iB,C,iC)$ is an orthonormal
basis for $\Re H$ and 
\[
\left<A,iA,B,C+\lambda iB\right>=W
\]
for some $\lambda\in\R$.

Then we choose
\begin{itemize}
\item
$L=\left<C\right>_\C$ if $|\lambda|>1$,
\item
$L=\left<B+iC\right>_\C$ if $0\le\lambda\le 1$ and
\item
$L=\left<B-iC\right>_\C$ if $-1\le\lambda <0$.
\end{itemize}

It is easy to check that in each case the angle is at least $\pi/4$.
\end{proof}

\begin{proposition}\label{prop-4}
Let $A$ be an abelian threefold (i.e.~an abelian variety of
dimension three).
Then there exists a smooth complex curve $C\subset A$
such that for every real $4$-dimensional subtorus $S\subset A$
and every point $a\in A$ there exists a point $p\in C$
where $C$ and $S(a)$ (the $S$-orbit in $A$ through $a$)
intersect transversally.
\end{proposition}
\begin{proof}
We have to consider all $4$-subtori. Since the set of all such tori
lacks good geometric properties, we instead consider the larger set
of all connected real Lie subgroups of real dimension $4$,
or, equivalently, the real Grassmann variety $M$ which
parametrizes all real vector subspaces of dimension $4$
of the Lie algebra $\Lie(A)\simeq\C^3$.
This is a real compact variety.

Now we fix an hermitian inner product on $\Lie(A)\simeq\C^3$
(e.g.~the standard one for $\C^3$ or the one corresponding to the
Riemann condition).
For each element $H\in M$ we define a closed neighbourhood $B_H$ as follows:
An element $H'$ belongs to $B_H$ iff for every vector $v'$ in $H'$
there is a vector $v\in H$ such that the angle between $v$ and $v'$
is at most $\pi/16$.
Due to compactness of $M$ there is a finite collection of elements
$H_i\in M$, $i\in I$ such that $M=\cup_{i\in I}B_{H_i}$.

Next we will choose a smooth complex curve $C_i\subset A$ for each
$i\in I$.
Fix an index $i$. Let $H=H_i$
and $B=B_{H_i}$.
Choose a complex line $L$ in $T_eA\simeq\Lie(A)$
such that the angle between $L$ and $H$ is at least $\frac{\pi}{4}$
(which is possible due to lemma~\ref{linalg}).
Then we choose a smooth complex curve $S=S_i$ through $e$ such that
for each $v\in T_eS$ there is a vector $v'\in L$ such that
the angle between $v$ and $v'$ is at most $\pi/16$
(lemma~\ref{bertini}).
By the definition of $B$, 
the angle between $T_eS$ and $H'$ is at least $\frac{\pi}{8}$ for
every $H'\in B$.

Now let $\pi:\C^3\to A$ denote the universal covering.
Let $F$ denote a fundamental region, i.e. a compact subset of $\C^3$
with $\pi(F)=A$. Let $W$ be an open neighbourhood
of $e$ in $S$ which is small enough such that the embedding of $W$ in
$A$ lifts to an embedding into $\C^3$, taking $e$ to $0$.
In addition, we require that $W$ is small enough such that
for every $w\in W$, $v\in T_wS\setminus\{0\}$
and $v'\in T_eS\setminus\{0\}$
the angle between $v$ and $v'$ is at most $\pi/16$.

For each $H'\in B$ we define 
\[
Z(H')=\{c+h:c\in W, h\in H'\}.
\]
We claim: There exists a number $\rho>0$ such that
$Z(H')$ contains the ball with radius $\rho$
and center $0$ for every $H'\in B$.
Indeed, assume the contrary. Then there are sequences $v^{(k)}\in\C^3$
and $H^{(k)}\in B$ such that:
\begin{enumerate}
\item
$H^{(k)}$ converges to an element $H''\in B$ (recall that $B$ is
compact),
\item
$\lim v^{(k)}=0$,
\item
$v^{(k)}\not\in Z(H^{(k)})$.
\end{enumerate}
But this would contradict the fact that $H''$ and $T_eW$ are
transversal.
Thus we can find such a number $\rho$.
Next, using compactness of $F$, we choose a finite set
$\Sigma\subset\C^3$
such that for every $x\in F$ there is an element
$s\in\Sigma$ with $||x-s||<\rho/3$.

Using this fact and lemma~\ref{union} we can find
a map $\xi:\Sigma\to\C^3$ such that:
\begin{itemize}
\item
$\cup_{s\in\Sigma}(\pi(\xi(s))+S)$ is smooth and
\item
$|\xi(s)|<\rho/3$ for all $s\in\Sigma$
\end{itemize}

Then we have constructed a smooth curve in $A$,
namely 
\[
S'=\cup_{s\in\Sigma}(\pi(\xi(s))+S)
\]
with the following property:

{\em (T) For every vector $u\in\C^3$ with $||u||<\rho/3$
and every real $4$-dimensional 
subtorus $H$ of $A$ with $\Lie(H)\in B$
every $H$-orbit in $A$ intersects $\pi(u)+S'$ in some point
transversally.
}

We found this curve $S'$ after fixing an element $i\in I$.
We can do the same for every element $i\in I$, obtaining a family
of curves $S'_i$ and a family of positive real numbers $\rho_i$.

Then by lemma~\ref{union} we can choose vectors $u_i$ such that
$||u_i||<\rho_i/3$ and such that
$C=\cup_{i\in I}\left( \pi(u_i)+S'_i\right)$ 
is a smooth curve.

By construction this curve has the property that it intersects each
translate of each real $4$-dimensional subtorus of $A$ in at least
one point transversally.
\end{proof}

\begin{proposition}\label{prop-ex}
There exists an abelian threefold $A$ with a complex curve $C$
such that the following property holds:

For every complex one-parameter subgroup $P$ of $A$ and every
point $a$ in $A$ there is a point $p$ in $C$ where
$C$ and the (real) closure of $P\cdot a$ intersect transversally.
\end{proposition}
\begin{proof}
We may choose $A$ such that every real three-dimensional real subtorus
is totally real (prop.~\ref{3-subtori}).
Then evidently real subtori of smaller dimension are
totally real as well. Now let $P$ be a complex one-parameter subgroup
of $A$. The closure of $P$ is again a subgroup, and therefore
in fact a real subtorus. This subtorus does not need to be
complex, but it can not be totally real, since it contains $P$.
Therefore for every complex one-parameter subgroup
$P$ of $A$ the real dimension of its closure is at least $4$.
Now it suffices to choose the curve $C$ according to thm.~\ref{prop-4}.
\end{proof}

\section{Brody curves and sets of rational points}

Conjecturally entire curves or Brody curves 
with values in projective varieties defined over 
some number field behave somewhat
analoguously to sets of rational points 
(admitting finite field extensions).

As we have seen, Brody curves and arbitrary entire curves
behave differerently. So which are the right analogue for rational
point sets? In our construction at one point we made a ``very
generic'' choice. For this reason it is not clear whether one can
find such an example which is defined over a number field.

If such an example can be defined over a number
field, it would suggest that complex-analytic concept corresponding
to infinite rational point sets are arbitrary entire curves
and not Brody curves: For every abelian variety $A$ defined
over a number field $k$ there is a finite field extension
$K/k$ such that $A(K)$ is Zariski dense. Then also $X(K)$ is
Zariski dense in $X$ for every projective manifold $X$ obtained 
from $A$ by blowing up something. 
Thus if our construction can be realized over a number field,
it would yield a projective variety defined over some number field $K$
such that every Brody curve is degenerate, but there is a Zariski
dense subset of $K$-rational points.

In any case, dense sets of rational points as well as dense entire
curves behave nicely under birational transformations while
our example shows that the behaviour of Brody curves may
change dramatically.

This suggests that the right complex-analytic analogue to infinite sets of rational
points should be arbitrary entire curves rather than Brody curves.

\end{document}